\documentclass[a4paper,12pt]{amsart}
\usepackage{subfigure}
\usepackage{graphicx}
\usepackage{rotating}
\usepackage{morefloats}
\usepackage[utf8]{inputenc}
\usepackage[T1]{fontenc}
\usepackage{amssymb}
\usepackage{amsmath}    
\def\R{I\kern -0,37 em R}
\def\P{I\kern -0,37 em P}
\def\Z{I\kern -0,37 em Z}

\begin{document}
\title[Lie and Cartan Theory, III]{On the Lie and Cartan Theory of Invariant Differential Systems, III}
\author{A. Kumpera}
\address[Antonio Kumpera]{Campinas State University, Campinas, SP, Brazil}
\email{antoniokumpera@gmail.com}

\date{March, 2017}

\keywords{Lie groupoids, Lie algebroids, Jordan-Hölder resolutions, differential systems, Pfaffian systems.}
\subjclass[2010]{Primary 53C05; Secondary 53C15, 53C17}

\maketitle

\begin{abstract}
It is presently our aim to undertake the discussion, of the Parts I and II, on the infinitesimal level and outline as well the transition from infinitesimal to finite, the main reason for this being, of course, the well known fact that arguments and calculation on the infinitesimal level are far simpler that those on the finite level.
\end{abstract}

\section{Introduction}
The Second and the Third Theorems of Lie do have drawbacks in the case of Lie groupoids, the whole context being considerably more complicated. We shall therefore begin by exhibiting a criterion, necessary and sufficient, for the appropriate integrability of a Lie algebroid, the same criterion holding for both Lie Theorems since, in ultimate analysis, the entire  difficulty arises when quotienting a foliation \textit{i.e.}, when considering the space of all the  \textit{leaves}.

\vspace{2 mm}
\noindent
Given a Lie groupoid $\mathcal{G}$ defined on a base manifold \textit{M} by a non-linear differential equation \textit{i.e.}, a groupoid contained in some jet space, there is associated to it a Lie algebroid defined on the same manifold by a Lie equation (\cite{Kumpera1972}). In this paper we intend to examine some aspects of the converse statement that is not always fulfilled. The main ingredients in our whole discussion are the so-called \textit{Cartan pseudo-groups} and \textit{Cartan groupoids} where everything works beautifully, Élie Cartan being the great master of the \textit{co-variant} universe.

\section{The differentiable quotientability of a foliation}
We adopt here, in discussing the integrability of a Lie algebroid, a geometric approach different from that found in \cite{Crainic2003} and begin by examining the quotientability of foliations. The following Lemma is a direct consequence of the construction of the \textit{holonomy groupoid} via the iterated juxtaposition of \textit{sliced} coordinate neighborhoods (\textit{c.f.},\cite{Chevalley1946}, \cite{Reeb1958},\cite{Hector1981},\cite{Mucuk1995},\cite{Aof1992},\cite{Mrcun2005}). We say that a groupoid is \textit{trivial} whenever, for any two unit elements, there is just a single element with such a source and such a target together with its inverse. We also recall that the equivalence relation canonically associated to a foliation is the one that renders equivalent any two points belonging to the same leaf and only these.

\vspace{2 mm}
\noindent
\newtheorem{alcatra}[LemmaCounter]{Lemma}
\begin{alcatra}
The holonomy groupoid of a foliation is trivial if and only if the equivalence relation associated to this foliation is a differentiable sub-manifold of the product manifold.
\end{alcatra}

\vspace{2 mm}
\noindent
On account of the Godement Theorem (\cite{Serre1965}) and its slight improvement brought by Kubarski (\cite{Kubarski1987}), we can state the following result where "\textit{differentiably quotientable}" means that the quotient of a foliation, \textit{modulo} its leaves, admits a differentiable ($e.g., C^{\infty}$) manifold structure for which the quotient map is a submersion. The quotient structure is then unique.

\vspace{2 mm}
\noindent
\newtheorem{alcatro}[TheoremCounter]{Theorem} \begin{alcatro}
In order that a foliation $\mathcal{F},$ defined on a differentiable manifold M, be differentiably quotientable, it is necessary and sufficient that the following two conditions hold:

\vspace{3 mm}
\hspace{1 mm}$(i)$ The holonomy groupoid associated to $\mathcal{F}$ is trivial and

\vspace{2 mm}
$(ii)$ the first (or, equivalently, the second) projection

\begin{equation*}
pr_1:\mathcal{R}~\longrightarrow~M/\sim
\end{equation*}

\vspace{2 mm}
\noindent
is a submersion, where $\mathcal{R}\subset M\times M$ denotes the equivalence relation associated to the foliation and $\sim$ the corresponding equivalence among the points in M. 
\end{alcatro}

\section{Jordan-Hölder resolutions for Lie algebroids}
Let us denote by $\mathcal{A}$ a Lie algebroid defined on the manifold \textit{M} (\cite{Kumpera2015}). By definition, $\mathcal{A}$ is a finitely generated free sub-sheaf of the sheaf $\underline{TM}$ of all the germs of vector fields, defined on $M,$ invariant under the Lie bracket that further satisfies the following restrictive \textit{differential} condition: It is the set of all the solutions of an integrable linear partial differential equation of a certain order $k,$ called a Lie equation, where after it will be said to be of order $k.$

\vspace{2 mm}
\noindent
We consider a Jordan-Hölder resolution
\begin{equation*}
\mathcal{A}=\mathcal{A}_0\supset\mathcal{A}_1\supset~\cdots~\supset\mathcal{A}_{\ell-1}\supset\mathcal{A}_{\ell}=0
\end{equation*}

\vspace{2 mm}
\noindent
for $\mathcal{A}$ where, by definitions, each $\mathcal{A}_{i+1}$ is a maximal sub-sheaf in ideals of $\mathcal{A}_i$ and, consequently, the quotient is simple. The integer $\ell$ is independent of the resolution and is called the Jordan-Hölder length of $\mathcal{A}.$ Occasionally, the initial term $\mathcal{A}$ will also be indicated, as above, by $\mathcal{A}_0.$ The resolution is called admissible when each term, belonging to it, is an integrable algebroid. The resulting sequence, at the groupoid level, will of course be a Jordan-Hölder resolution for the groupoid $\mathcal{G}$ obtained by integrating $\mathcal{A}~.$  
\vspace{2 mm}
\noindent
Given a Lie groupoid $\mathcal{G},$ the foliation defined by its $\alpha-$fibres \textit{i.e.}, the collection of all the connected components of every $\alpha-$fibre is a differentiably quotientable foliation on the underlying manifold to the Lie groupoid $\mathcal{G}$ and, in particular, when the $\alpha-$fibres are connected, the quotient manifold is equal (diffeomorphic) to the sub-manifold of units of $\mathcal{G}.$ The same comment can be paraphrased with respect to the $\beta-$fibres.

\vspace{2 mm}
\noindent
Conversely, let $\mathcal{A}$ be an integrable Lie algebroid and let us denote by $\mathcal{G}_0$ the Lie groupoid obtained by extending globally, according to Pradines (\cite{Almeida1980},\cite{Pradines1966}), the "\textit{piece of differentiable groupoid}" (\textit{morceau de groupoïde différentiable}) obtained by integrating, for example, the distribution determined by the \textit{right invariant} vector fields and by considering just the leaves that contain some element of \textit{M}. Taking a Jordan-Hölder resolution of the algebra $\mathcal{A}~,$ we proceed inasmuch with each sub-algebroids, taking only the leaves that contain some element of \textit{M}, and obtain a differentiable (Lie) sub-groupoid with connected $\alpha-$fibres and with the same space of units. Let us state this property in a more suitable fashion.

\vspace{2 mm}
\noindent
\newtheorem{alcacuz}[LemmaCounter]{Lemma}
\begin{alcacuz}
When $\mathcal{A}$ is an integrable algebroid, then each term belonging to a Jordan-Hölder resolution of $\mathcal{A}$ is also integrable.
\end{alcacuz}

\vspace{2 mm}
\noindent
\textbf{Remark.} The above property also holds inasmuch for the terms of any composition series of $\mathcal{A}~.$

\section{The Lie algebroids associated to a differential system.}
We now consider a differential system $\mathcal{S}$ of some order \textit{k} , considered as a (locally trivial) sub-bundle of the $k-$th order jet space $J_k\pi,$ where $\pi:P~\longrightarrow~M$ is a fibration, and denote by $\mathcal{C}_k$ the canonical contact structure at that same order \textit{k} (\cite{Kumpera2016}). It is well known that a local vector field, defined on the space of all the invertible $k-$jets, is locally the prolongation of a vector field defined in the base space if and only if it leaves invariant the $k-$th order canonical contact structure. This being so, we can associate to the differential system $\mathcal{S},$ given above, the sheaf of all the germs of vector fields defined on the base space whose prolongations are tangent to the sub-manifold $\mathcal{S}.$ Actually, it will suffice that the prolongued vector fields be tangent, to order 1, to $\mathcal{S}$ at the points in consideration since the first order tangency will imply the $0-$th order tangency of the bracket. On account of the above mentioned relationship between prolongued vector fields and those leaving invariant the contact structure, we can restate all the above considerations exclusively in terms of the germs of vector fields defined on the space of $k-$jets.

\vspace{2 mm}
\noindent
The algebroid $\mathcal{A}~,$ defined above, is not necessarily integrable but we can determine a \textit{maximal} integrable sub-algebroid $\mathcal{A}_0~.$ In fact, given any germ of vector field tangent, after prolongation, to $\mathcal{S}$ at the corresponding points, it generates, according to Pradines, a piece of differentiable groupoid that can be globalized onto a Lie groupoid leaving invariant the differential system (the bracket of a vector field with itself is null). We can therefore take the sum of all the Lie sub-algebroids that generate global Lie sub-groupoids and obtain the maximal integrable sub-algebroid $\mathcal{A}_0$ associated to the differential system $\mathcal{S}.$ Let us denote by $\mathcal{G}_0$ a piece of differentiable sub-groupoid that it generates and by $\mathcal{G}$ the corresponding global Lie sub-groupoid resulting 
thereafter. It is the largest Lie groupoid that leaves invariant the differential system at that order.

\vspace{3 mm}
\noindent
Let us now recall that the purpose, for considering Jordan-Hölder sequences, resides in their relevance with respect to the integration of differential systems via the Lie-Cartan algorithm consisting in successive reductions, of the given system, to auxiliary systems invariant under simple Lie groupoids (\textit{c.f.}, \cite{Kumpera2016}). However, in order to do so, we are forced initially to integrate the terms of the given infinitesimal series of the Lie algebroid and only, thereafter, apply the Lie-Cartan method with the help of the finite resolution obtained upon integration. There is, nevertheless, an alternative approach that seems somehow more promising. Instead of considering directly the given differential system, we take its linearization, apply the Lie-Cartan method to this linear differential system with the help of the infinitesimal series and finally obtain the solutions of the given system by integrating those of the linear system. It is somehow a \textit{Tom and Jerry} mischievous game.

\vspace{2 mm}
\noindent
Given a differential system $\mathcal{S}$, we consider the \textit{tangent fibration}
\begin{equation*}
\pi\circ p:TP~\longrightarrow~P~\longrightarrow~M~,
\end{equation*}
where $p:TP~\longrightarrow~P$ is the natural projection. Furthermore, $J_kTP=J_k(\pi\circ p)$ projects onto $J_kP=J_k\pi$ via $j_kp$ , each element $X\in J_{k+1}TP$ projects onto a tangent vector of $J_kP$ at the point $j_kp(X)$ ($J_{k+1}\subset J_1J_k$) and the fibration
\begin{equation*}
\beta\circ j_{k+1}p:J_{k+1}TP~\longrightarrow~P
\end{equation*}

\vspace{2 mm}
\noindent
is a vector bundle. We can now associate, to the $k-$th order equation $\mathcal{S},$ the linear equation of order $k+1,$ denoted by $\mathcal{LS}$ and composed by the set of all the elements in $J_{k+1}TP$ whose images under $\rho_{k,k_1}\circ j_{k+1}p$ are tangent to $\mathcal{S}.$

\vspace{2 mm}
\noindent
Any local vector field defined on $P$ prolongs to TP hence also to $J_{k+1}TP$ where this latter jet space has the above described meaning, not to be misunderstood as being the space of $(k+1)-$st jets of vector fields on the manifold $P.$ There after, we can determine the Lie algebroid defined on $P$ and consisting of all the germs of vector fields that, after prolongation, leave invariant the linear equation $\mathcal{LS}$ \textit{i.e.}, are tangent to the equation or, equivalently, generate local $1-$parameter groups that preserve the equation. Taking a Jordan-Hölder resolution for the above Lie algebroid, the Lie-Cartan algorithm can now be applied inasmuch as it is described in \cite{Kumpera2016} since the $\alpha-$fibres (and the $\beta-$fibres) of a Lie groupoid are no other that the integral manifolds of its Lie algebroid considered as its set of right invariant vector fields (resp. left invariant). Once obtained the solutions of the linear equation $\mathcal{LS}$, we obtain those of the general equation $\mathcal{S}$ simply projecting these solutions by means of  $p$ \textit{i.e.}, composing them with $p.$

\vspace{3 mm}
We now briefly recall the construction of a piece of differentiable groupoid, corresponding to a given Lie algebroid.

\vspace{2 mm}
\noindent
Let us consider a Lie algebroid $\mathcal{L}$ defined on the base manifold $M$ and let us  denote by $\mathcal{H}$ the set of all the homotopy classes of integral curves of $\mathcal{L}$ that initiate and terminate at the same points.\footnote{Integral curves of representatives of the elements belonging to $\mathcal{L}.$} If $\gamma(t)$ is an integral curve of $\xi$ then, for any constant $c,~\gamma(ct)$ is an integral curve of $c\xi$ hence we could restrict our considerations to the sole integral curves defined on the unit interval though this does not seem to contribute with any mayor simplifications. Given such a homotopy class $[\gamma],$ $\gamma:[0,t_0]~\longrightarrow~M,$ we shall say that $\gamma(0)$ is the \textit{source} of $[\gamma],$ denoted by $\alpha([\gamma]),$ and that $\gamma(t_0)$ is the target, denoted by $\beta([\gamma])$.Next, we consider chains (sequences) of elements of $\mathcal{H}$ namely, objects of the form 
\begin{equation*}
([\gamma_p],[\gamma_{p-1}],~\cdots~,[\gamma_1])~, 
\end{equation*}

\vspace{2 mm}
\noindent
where $\beta([\gamma_i])=\alpha([\gamma_{i+1}]).$ The set of all such chains, for arbitrary values of the integer $p,$ then becomes a groupoid $\mathcal{G}$ under the obvious composition of chains, the inverse of an element being that element determined by the same path with the opposite (contrary sense) parametrization or a composite of such elements as soon as we remark that the source of the inverse is the target of the element and inasmuch for the target. Finally, we identify two chains that have the same source and the same  target. This identification is an equivalence relation compatible with the groupoid structure of the chains hence we obtain again a groupoid whose $\alpha-$fibres are represented by sub-sets of $M$ that, however should not be identified with such sub-sets, since the same point of $M$ can represent distinct points in different $\alpha-$fibres of the above groupoid. The property of $\mathcal{L}$ being free ensure that the distribution generated by all the above vector fields (and always in the vicinity of $M$ represented now by the constant curves) is regular and integrable where after we can choose, at each point of $M,$ a coordinate system adapted to this foliation and, in particular, a \textit{cubic} neighborhood in the sense of Chevalley (\cite{Chevalley1946}). All such cubic neighborhoods then provide the desired differentiable structure of $\mathcal{G}$ in the vicinity of $M,$ this structure admitting eventually some extension beyond these cubic neighborhoods. 

\section{Lie Equations.}
Let $\pi:E~\longrightarrow~M$ be a \textit{fibration i.e.}, the mapping $\pi$ is a differentiable  surjective submersion. We denote by $J_kE$ or by $J_k\pi,$ when it is desirable to put in evidence the projection $\pi,$ the manifold of $k-$jets of local sections of $\pi$ (\textit{Ehresmann, \cite{Ehresmann1967}}), by $\alpha$ and $\beta$ the source and target projections respectively and by 

\vspace{3 mm}
\noindent
\begin{equation*}
\alpha\times\beta:J_kE~\longrightarrow~M\times E
\end{equation*}

\vspace{2 mm}
\noindent
their fibered product.

\vspace{2 mm}
\noindent
\newtheorem{alcaxofra}[LemmaCounter]{Lemma}
\begin{alcaxofra}
When E is a vector bundle over M then so is $J_kE$ with respect to the projection $\alpha$ onto M.
\end{alcaxofra}

\vspace{2 mm}
\noindent
The proof is obvious since it suffices to define the sum of two $k-$jets as the $k-$jet of the sum and inasmuch for the scalar multiplication.

\vspace{5 mm}
\noindent
We now consider the tangent bundle $TM$ fibrating over \textbf{M} and denote by $j_k\sigma$ the section of $J_kTM$ obtained by taking all the $k-$jets of the local section $\sigma$ of $TM~\longrightarrow~M.$ Since the definition of the bracket of two vector fields requires the first order derivatives of their components, we conclude that, given the $k-$jets $j_k\xi(x)$ and $j_k\eta(x),$ the $(k-1)-$th jet $j_{k-1}[\xi,\eta](x)$ of the bracket is well defined and will be denoted by $[j_k\xi(x),j_k\eta(x)].$ 

\vspace{2 mm}
\noindent
\newtheorem{alcato}[DefinitionCounter]{Definition}
\begin{alcato}
A (locally trivial) vector sub-bundle $\mathcal{L}$ of $J_kE$ is called a Linear Lie Equation of order k when the bracket of any two $k-jets,$ contained in the equation, belongs to the $(k-1)-$st projection $\rho_{k-1,k}(\mathcal{L})$ of $\mathcal{L}.$  
\end{alcato}

\vspace{2 mm}
\noindent
The most outstanding and, for that matter, the most useful property of a Lie equation states that the bracket of any two local vector fields that are solutions of the equation is also a solution. When the equation is integrable then, conversely, the above property characterizes its Lie nature. An analytic Lie equation is always integrable.

\vspace{2 mm}
\noindent
The notation $\rho_{h,k}:J_kTM~\longrightarrow~J_hTM$ is of course standard. We shall now modify appropriately the Lie bracket of vector fields so as to be able not only to define but as well to manipulate more adequately a Lie Equation that will be re-defined as an equation invariant under the modified bracket.

\vspace{2 mm}
\noindent
Let us consider two differentiable manifolds $M$ and $N$ of the same dimensions and denote by $\Pi_k(M,N)$ the manifold of all the invertible $k-$jets originating in $M$ and terminating in $N.$ When $M=N,$ we then just write $\Pi_kM.$ The fibered product $\alpha\times\beta$ is then a surmersion onto $M\times N.$ We could, of course, define $\Pi_kM$ in terms of the $k-$jets of local sections of a fibration but our definition seems more appropriate and easier to handle. The set $\Pi_kM$ is a Lie groupoid and $\Pi_k(M,N)$ is a left/right homogeneous space of the corresponding Lie groupoids. We also recall that the Lie algebroid associated to $\Pi_kM$ is equal (isomorphic) to the sheaf of germs of local sections of $J_kTM.$

\vspace{2 mm}
\noindent
We denote by $M^2=M\times M$ the product manifold, $\mathcal{F}_{M^2}$ its structure sheaf (germs of differentiable functions), $\mathcal{O}_{M^2}$ the sub-sheaf, in ideals, defining the diagonal (germs that vanish on the diagonal) and $\Delta$ the diagonal. We also call \textit{diagonal} the vector fields on $M^2$ that are tangent to the diagonal. The first basic remark states that $\mathcal{O}_{M^2}$ is invariant under the Lie derivative with respect to a diagonal vector field \textit{i.e.}, $\theta(\zeta)\mathcal{O}_{M^2}\subset\mathcal{O}_{M^2},$ for any diagonal vector field $\zeta.$ In order to understand our further reasoning, the above property should be confronted with the behaviour, at a point $x\in M,$ of the ideal of those functions that vanish at $x.$ The Lie derivative of such a function with respect to a vector field that vanishes at the point $x$ does not necessarily vanish at that point. We next define a \textit{diagonal} module structure on the Lie algebra $\Xi(M^2)$ of all the vector fields, defined on the product manifold, by setting
\begin{equation*}
f\Delta\zeta=g\zeta_H+f\zeta_V
\end{equation*}

\vspace{2 mm}
\noindent
On account of the previously mentioned \textit{basic remark}, we infer that the $k-$jet, at a point of the diagonal, of the bracket of two diagonal vector fields only depends and is entirely determined by the $k-$jets of these two vector fields not depending any more of the $(k+1)-$jets. We finally select a convenient $\Delta-$sub-module of $\Xi(M^2)$ namely, the set $\Theta$ of all the vector fields that are $\pi_1-$projectable onto $M$. In other terms, we just take those diagonal vector fields $\zeta$ that verify
\begin{equation*}
T_{(x,y)}\pi_1(\zeta_H)=T_{(x,y')}\pi_1(\zeta_H)     
\end{equation*}

\vspace{2 mm}
\noindent
at arbitrary points $x,y,y'.$ We next observe that this $\Delta-$sub-module, denoted by $\Theta,$ is also a sub-algebra with respect to the Lie algebra structure of the set $\chi$ of all the vector field on the product space and the above remark shows that the Lie bracket of vector fields induces a (fibre-wise) Lie bracket on the $k-$jets of the elements of $\Theta,$ along the diagonal. However, the projection $\pi_1$ extend, along the diagonal, to a surjective projection of the space of $k-$jets of the elements of $\Theta,$ along the diagonal, onto the space of $k-$jets of vector fields on $M$ \textit{i.e.}, onto the space $J_kTM.$ Since the bracket operation is preserved by projection, we finally conclude that the jet space $J_kTM$ is provided with a modified Lie bracket $\textbf{[~,~]}$ that satisfies the following property:

\vspace{2 mm}
\newtheorem{grange}[PropositionCounter]{Proposition}
\begin{grange}
The equation $\mathcal{L}\subset J_kTM$ is a Lie equation if and only if it is invariant under the modified Lie bracket.
\end{grange}

\vspace{2 mm}
\noindent
The proof is a direct consequence of the fact that the horizontal components of the elements of $\Theta$ are simply obtained by \textit{horizontally propagating}, to $M^2,$ the vector fields on $M.$ For further details, we refer the reader to \cite{Kumpera1972}.

\section{The Contact Ambient}
In order to fix the terminology and avoid any conflicts, let us state the following. A \textit{pseudo-group} is a set of local transformations, a \textit{pseudo-algebra} is a set of local vector fields, a \textit{groupoid} is a set with a partially defined multiplication operation, an \textit{algebroid} is a sheaf of germs of vector fields and finally a \textit{Lie equation} is a linear partial differential equation whose local sections can be composed via a Lie bracket. We provided, in the previous section, a brief overview on linear Lie equations though a much more thorough and complete study of these can be found in the references, especially in \cite{Kumpera1972}. At present, we shall only assume a differentiability class high enough so as to enable us to perform the necessary calculations.

\vspace{2 mm}
\noindent
We first consider a Lie groupoid $\mathcal{G}$ of order $k+1$, defined in the Jet space $\Pi_{k+1}M$ of invertible $(k+1)-$jets on a manifold \textit{M} of dimension equal to \textit{n}, and look for the appropriate spaces where it could be conveniently placed. Since every $(k+1)-$jet $X\in\Pi_{k+1}M$ can be lifted (prolonged) to a $1-$jet in $\Pi_kM$ originating at any jet whose source is equal to the source of \textit{X}, the groupoid $\mathcal{G}_{k+1}~,$ of order $k+1,$ becomes a groupoid $\mathcal{G}_{1,k}~,$ of order 1, defined on $\Pi_kM$ and isomorphic to the original one. Inasmuch, given a pseudo-group of order $k+1,$ we can extend its elements (local transformations of \textit{M}) to $\Pi_kM$ and thus obtain a first order pseudo-group on this manifold. Moreover, its first order groupoid is precisely the groupoid obtained from its $(k+1)-$st order groupoid. Let us denote by $\mathcal{C}$ or, when necessary, by $\mathcal{C}_k$ the canonical contact structure defined, on $\Pi_kM~,$ by the first degree contact forms. It is well known that the prolongation of every local transformation coming from any base space is a local automorphism of this contact system. Then the canonical contact structure $\mathcal{P}_k$ can be restricted to $\mathcal{G}_k$ and all the geometry inherent to the groupoid can be investigated by means of this restricted Pfaffian system. In particular, the extension of every jet $X\in\mathcal{G}_k~,$ to order $k-1~,$ defines a first order Lie groupoid canonically isomorphic to the given groupoid.

\vspace{2 mm}
Let us also observe that everything that was said for pseudo-groups and groupoids transcribes \textit{ispis litteris} for pseudo-algebras (or algebroids) and \textit{Lie equations} (\cite{Kumpera1972}).

\vspace{2 mm}
We next consider a Lie pseudo-group $\Gamma~,$ of order \textit{k} , defined on the same manifold \textit{M}, and denote by $\Gamma_k$ the corresponding $k-$th order (Lie) groupoid. Then $\Gamma$ prolongs to a first order  pseudo-group $\tilde{\Gamma}$ defined on the manifold $\Pi_kM$ and its restriction to the sub-manifold $J_k\Gamma$ composed by all the $k-$jets of its elements (local transformations of \textit{M}) becomes a \textit{Cartan} pseudo-group (\textit{cf}.\cite{Kumpera1963}).

\vspace{2 mm}
Let us finally state the following result which is, in fact, rather obvious since the generated \textit{finite object} is simply obtained by means of the same \textit{invariance forms} as those defining the algebroid. The proof is only required to assure that the finite pseudo-group, defined as above, is large enough so as to display the entire algebroid. Cartan also gave two proofs where, however, he shows much deeper results since he starts with constants that should eventually become those pertaining to the structure equations relative to Pfaffian forms to be constructed (\cite{Cartan1904},\cite{Cartan1905},\cite{Cartan1937}).

\vspace{2 mm}
\newtheorem{jhi}[TheoremCounter]{Theorem}
\begin{jhi}
To any Cartan pseudo-algebra (resp. algebroid) corresponds a Cartan pseudo-group (resp. a groupoid) with a pseudo-algebra (resp. an algebroid) equal to the given one. 
\end{jhi}

\vspace{2 mm}
This being so, we can associate to any Lie algebroid a Lie groupoid inasmuch as it happens in Lie group theory. The only drawback resides in the fact that this Lie groupoid or, for that matter, the corresponding Lie pseudo-group do not sit on the desired manifold \textit{M}. Our task is now to find suitable conditions, on the algebroid or on the Lie equation, so as to be able to project, onto \textit{M}, the above pseudo-group as well as the groupoid. A first obvious remark is the following:

\vspace{2 mm}
\noindent
\textit{The local transformations obtained by the  integration of the Lie equation should be contact transformations since then, they would project locally and the local projections of a given transformation can, of course, be glued together.}

\section{Lie equations compatible with the contact structures.}
Needless to say what we mean by the \textit{compatibility}. Since any local transformation obtained by the integration of a solution of a Lie equation $\mathcal{L}$ is a local automorphism of $\mathcal{L}~,$ the condition transcribes by $Aut\mathcal{L}\subset Aut\mathcal{C}~,$ where $\mathcal{C}$ denotes the above considered contact structure. We start by examining a little closer the bracket involving Lie equations.

\vspace{2 mm}
\noindent
Lie equations are not algebroids but are the defining equations of those algebroids composed by germs of sections of a tangent bundle $TM~\longrightarrow~M.$ In trying to integrate an algebroid, one is forced to determine the non-linear Lie equation associated to the above linear Lie equation (defining the algebroid) and then see what can be done with the integration of this non-linear Lie equation. Its solutions provide what Pradines calls a kernel (noyau) or piece (morceau) of a Lie groupoid and one must search the eventual globalization of this kernel. On account of the previous discussion, it suffices to consider the groupoid $\mathcal{G},$ of section 4, constructed with the help of chains and verify whether the distribution defined by the algebroid, whose integrals are either the $\alpha$ or the $\beta$ fibres, is differentiably quotientable. The Theorem 1 provides a reasonably geometric criterion.

\vspace{2 mm}
\noindent
Finite dimensional real Lie algebras are all integrable but, in the case of Lie algebroids, most of them are not. Let us give a simple example where we consider one of the differential systems exhibited in \cite{Kumpera2016}, pg.15 and concerning the case (\textit{b}) where the second order partial derivative $\frac{\partial^2\textbf{S}}{\partial t^2}$ does not vanish identically. Let us assume that the full (transitive) invariance Lie algebroid is integrable. Then the resulting Lie groupoid $\mathcal{G}$ is also transitive and has an equivalent (isomorphic) Cartan representative of some order, composed of contact transformations jets. However, such a Cartan groupoid is defined by a finite set of Pfaffian forms satisfying a certain structure equation hence, the given algebroid has an isolated representative and consequently $\mathcal{G}$ could not belong to a continuous family.  

\vspace{2 mm}
\begin{figure}[htbp!]
\centering
\includegraphics[scale=3.0]{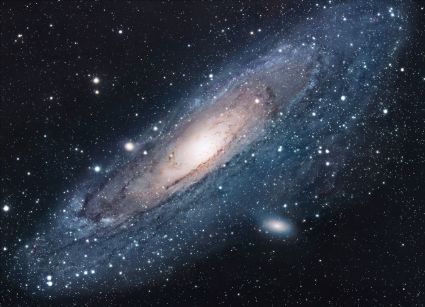}
\caption{The Universe}
\label{fig:univerise}
\end{figure}

\bibliographystyle{plain}
\bibliography{references}
\end{document}